\documentclass[11pt]{article}
\usepackage{amsmath}
\usepackage{amssymb}

\renewcommand{\Bbb}{\mathbb}
\newcommand{\B}{{\Bbb  B}}
\newcommand{\C}{{\Bbb  C}}
\newcommand{\D}{{\Bbb D}}
\newcommand{\Cn}{{\Bbb  C\sp n}}

\newcommand{\R}{{\Bbb  R}}

\newcommand{\Z}{{\Bbb  Z}}

\newcommand{\T}{{\Bbb  T}}

\newcommand{\F}{{\cal F}}

\newcommand{\I}{{\cal I}}

\renewcommand{\O}{{\cal O}}

\newcommand{\set }[1]{{\{#1\}}}

\newcommand{\PSH}{{\operatorname{PSH}}}

\renewcommand{\dim}{{\operatorname{dim}}}
\newcommand{\codim}{{\operatorname{codim}}}
\newcommand{\Sing}{{\operatorname{Sing}}}

\numberwithin{equation}{section}

\newtheorem{theorem+}           {Theorem}      [section]
\newtheorem{definition+}  [theorem+]  {Definition}
\newtheorem{lemma+}  [theorem+]  {Lemma}
\newtheorem{corollary+}  [theorem+]  {Corollary}
\newtheorem{proposition+}  [theorem+]  {Proposition}
\newtheorem{example+}  [theorem+]  {Example}

\newenvironment{theorem}{\begin{theorem+}\sl}{\end{theorem+}\rm}
\newenvironment{definition}{\begin{definition+}\rm}{\end{definition+}\rm}
\newenvironment{lemma}{\begin{lemma+}\sl}{\end{lemma+}\rm}
\newenvironment{corollary}{\begin{corollary+}\sl}{\end{corollary+}\rm}
\newenvironment{proposition}{\begin{proposition+}\sl}{\end{proposition+}\rm}
\newenvironment{example}{\begin{example+}\rm}{\end{example+}\rm}
\newenvironment{proof}{\medbreak\noindent{\it Proof:}\rm}{\hfill$\blacksquare$\rm}
\newenvironment{prooftx}[1]{\medbreak\noindent{\it #1:}\rm}{\hfill$\blacksquare$\rm}

\title{\Large \bf GREEN FUNCTIONS WITH SINGULARITIES
ALONG COMPLEX SPACES}
\author{\large\bf  Alexander Rashkovskii and  Ragnar Sigurdsson }

\date{August 30,  2004}
\begin{document}
\maketitle

\begin{abstract}  \noindent
We study properties of a Green function $G_{A}$ with singularities
along a complex subspace $A$ of a complex manifold $X$. It is
defined as the largest negative plurisubharmonic function $u$
satisfying locally $u\leq \log|\psi|+C$, where
$\psi=(\psi_1,\dots,\psi_m)$, $\psi_1,\dots,\psi_m$ are local
generators for the ideal sheaf $\I_A$ of $A$, and $C$ is a
constant depending on the function $u$ and the generators.  A
motivation for this study is to estimate global bounded functions
from the sheaf $\I_A$ and thus proving a ``Schwarz Lemma'' for
$\I_A$.
\end{abstract}

\section{Introduction}\label{sec:intro}

If $\varphi$ is a bounded holomorphic function on a complex
manifold $X$, then it is a natural problem to estimate $|\varphi|$
given  some information on the location of the  zeros of $\varphi$
and  their multiplicities. If  $|\varphi|\leq 1$ and the only
given information is that $\varphi(a)=0$ for a single point $a$,
then
$$
\log|\varphi|\leq G_{X,a}=G_X(\cdot,a),
$$
where $G_{X,a}$ is the {\it pluricomplex Green function with
logarithmic pole at } $a$.  It is defined as the supremum over the
class $\F_{X,a}$ of all negative plurisubharmonic functions $u$
such that $u\leq \log|{\zeta}|+C$ near $a$, where ${\zeta}$ are
local coordinates near $a$ with $\zeta(a)=0$ and $C$ is a positive
constant depending on $u$ and $\zeta$.  The function $G_{X,a}$ was
introduced and studied by several authors \cite{Lem}, \cite{Zah},
\cite{Kl1}, \cite{PSh}, \cite{D0}, see also \cite{Kl}, \cite{D1}.

A generalization is to take $A=(|A|,\set{m_a}_{a\in |A|})$, where
$|A|$ is a finite subset of $X$,  $m_a$ is a positive real number
for every $a\in |A|$, and assume that $\varphi$ has a zero of
multiplicity at least $m_a$ at every point $a$ in $|A|$.  Then
$$
\log|\varphi|\leq G_{A},
$$
where $G_{A}$ is the {\it Green function with several weighted
logarithmic poles}. It is defined as the supremum over the class
of all negative plurisubharmonic functions $u$ on $X$ satisfying
$u\leq m_a\log|\zeta_a|+C$ for every $a$ in $|A|$, where $\zeta_a$
are local coordinates near $a$ with $\zeta_a(a)=0$ and $C$ is a
positive constant depending on $\zeta_a$ and $u$. The function
$G_{A}$ was first introduced by Zaharyuta \cite{Zah} and
independently by Lelong \cite{Le}.

The notion of multiplicity of a zero of an analytic function has a
natural generalization as a Lelong number of a plurisubharmonic
function.  If $u$ is plurisubharmonic in some neighbourhood of the
origin $0$ in $\C^n$, then the {\it Lelong number} $\nu_u(0)$ of
$u$ at  $0$ can be defined as
$${\nu}_u(0)=\lim_{r\to 0}
\frac{\sup\,\{u(x);\: |x|\leq r\}}{\log r}$$ and if $u$ is
plurisubharmonic on a manifold $X$ then the Lelong number
$\nu_u(a)$ of $u$ at $a\in X$ is defined as $\nu_u(a)=\nu_{u\circ
\zeta^{-1}}(0)$, where $\zeta$ are local coordinates near $a$ with
$\zeta(a)=0$. It is clear that this definition is independent of
the choice of local coordinates and that $\nu_u(a)$ equals the
multiplicity of $a$ as a zero of the holomorphic function
$\varphi$ in the case $u=\log|\varphi|$. Note that the
pluricomplex Green function $G_{X,a}$ with logarithmic pole at $a$
can be equivalently defined as the upper envelope of all negative
plurisubharmonic functions $u$ on $X$ satisfying $\nu_u(a)\ge 1$
and, similarly, $\nu_u(a)\ge m_a$ for the Green functions with
several weighted logarithmic poles.

For any non-negative function $\alpha$ on $X$, L\'arusson and
Sigurdsson \cite{LarSig1}, \cite{LarSig2} introduced the Green
function $\tilde G_\alpha$ as the supremum over the class $\tilde
\F_\alpha$ of all negative plurisubharmonic functions  with
$\nu_u\geq \alpha$. It is clear that if $\varphi$ is holomorphic
on $X$, $|\varphi|\leq 1$,  and every zero $a$ of $\varphi$ has
multiplicity at least
 $\alpha(a)$, then
$$
\log|\varphi|\leq \tilde G_\alpha.
$$
In this context it is necessary to note that we assume that the
manifold $X$ is connected, we take the constant function $-\infty$
as plurisubharmonic, and set $\nu_{-\infty}=+\infty$.  Hence
$-\infty\in \tilde \F_\alpha$ for every $\alpha$.  By
\cite{LarSig1}, Prop.~5.1, $\tilde G_\alpha\in \tilde \F_\alpha$.

In the special case when $X$ is the unit disc $\D$ in $\C$, we
have
$$
\tilde G_\alpha(z)=\sum_{w\in \D} \alpha(w) G_\D(z,w), \qquad z\in
\D,
$$
where $G_\D$ is the Green function for the unit disc,
$$
G_\D(z,w)=\log\bigg|\dfrac{z-w}{1-\bar w z}\bigg|, \qquad z,w\in
\D.
$$

If $X$ and $Y$ are complex manifolds, $\alpha$ is a non-negative
function on $X$, and $\Phi:Y\to X$ is a holomorphic map, then the
pullback $\Phi^*u=u\circ \Phi$ satisfies $\nu_{\Phi^*u}\geq
\Phi^*\nu_u$, so  $\Phi^*u\in \tilde\F_{\Phi^*\alpha}$ for every
$u\in \tilde\F_\alpha$.  This implies $\Phi^*\tilde G_\alpha\leq
\tilde G_{\Phi^*\alpha}$, i.e., $\tilde G_\alpha(x)\leq \tilde
G_{\Phi^*\alpha}(y)$ if $x=\Phi(y)$, and in particular
$$\tilde G_\alpha(x)\leq \tilde G_{f^*\alpha}(0)=\sum_{w\in \D}
f^*\alpha(w)\log|w|, \qquad f\in \O(\D,X), \quad  f(0)=x.
$$  One of the main results of \cite{LarSig1} and
\cite{LarSig3} is that for every manifold $X$ and every
non-negative function $\alpha$ we have the formula
\begin{equation}
\tilde G_\alpha(x)=\inf\set{\tilde G_{f^*\alpha}(0) \,;\, f\in
\O(\overline \D,X), f(0)=x}, \qquad x\in X. \label{eq:disc}
\end{equation}
Here $\O(\D,X)$ is the family of all analytic discs in $X$ and
$\O(\overline \D,X)$ is the subclass of closed analytic discs,
i.e., maps from $\D$ to $X$ that can be extended to holomorphic
maps in some neighbourhood of the closed disc $\overline \D$.
Results of this kind originate in Poletsky's theory of analytic
disc functionals, started in \cite{PSh} and \cite{P0}.

A natural way  of describing the zero set of a holomorphic
function $\varphi$ is to state that its germs $(\varphi)_x$ are in
the stalk $\I_{A,x}$ of a prescribed coherent ideal sheaf
$\I_A=(\I_{A,x})_{x\in X}$  of a closed complex subspace $A$ of
$X$.  Then, if $\psi_1,\dots,\psi_m$ are local generators of
$\I_A$ near the point $a$,  the function $\varphi$ can be
represented as $\varphi=\varphi_1\psi_1+\cdots+\varphi_m\psi_m$
near $a$, which implies that $\log|\varphi|\leq \log|\psi|+C$ near
$a$, where $\psi=(\psi_1,\dots,\psi_m)$, $|\cdot|$ is the
euclidean norm, and $C$ is a constant depending on $\varphi$ and
the generators.

We define $\F_{A}$ as the class of all negative plurisubharmonic
functions $u$ in $X$ satisfying $u\leq \log|\psi|+O(1)$ locally in
$X$, and we define the function $G_{A}$, the {\it pluricomplex
Green function with singularities along $A$}, as the supremum over
the class $\F_{A}$.

It follows from the definition of $G_A$ that if $A'$ is the
restriction of $A$ to a domain $X'$ in $X$, $\varphi$ is
holomorphic function on $X$, and $(\varphi)_x\in \I_{A,x}$ for all
$x\in X'$, then
$$|\varphi|\le e^{G_{A'}(x)}\sup_{X'}|\varphi|,\quad x\in X',$$
which is a variant of the Schwarz lemma for the ideal sheaves.

In order to relate $G_{A}$ to the Green functions $\tilde
G_\alpha$ above, we  define the function $\tilde\nu_A$ on $X$ by
$\tilde\nu_A(x)=\nu_{\log|\psi|}(x)$ if
$\psi=(\psi_1,\dots,\psi_m)$ are local generators for  $\I_A$ in
some neigbourhood of $x$.  It is easy to see that $\tilde\nu_A(x)$
is independent of the choice of the generators (actually, it
equals the minimal multiplicity of the functions from $\I_{A,x}$
at $x$), so $\tilde\nu_A$ is a well defined function on $X$ and
$\nu_u\geq \tilde\nu_A$ for all $u\in \F_{A}$. Hence, with
$\tilde\nu_A$ in the role of $\alpha$ above, we have
$\F_{A}\subseteq \tilde \F_{\tilde\nu_A}$ which implies
$$
G_{A}\leq \tilde G_{\tilde\nu_A}.
$$
In general, $G_{A}\neq\tilde G_{\tilde\nu_A}$ as seen from the
example where $X=\D^2$ and $\I_A$ has the global generators
$\psi=(\psi_1,\psi_2)$ with $\psi_1(z)=z_1^2$ and $\psi_2(z)=z_2$.
Then  $G_{A}(z)=\max\set{2\log|z_1|,\log|z_2|}$ and $\tilde
G_{\tilde\nu_A}(z)=\max\set{\log|z_1|,\log|z_2|}$ for
$z=(z_1,z_2)\in \D^2$. If, on the other hand, $A$ is an effective
divisor generated by the function $\psi$ in an open subset $U$ of
$X$, then by \cite{LarSig2}, Prop.~3.2, the function $\tilde
G_{\tilde\nu_A}-\log|\psi|$ on $U\setminus |A|$ can be extended to
a plurisubharmonic function on $U$. This implies that
$G_{A}=\tilde G_{\tilde\nu_A}$ for effective divisors $A$.

Now to the content of the paper.  In Section~2 we present the main
results, which are proved in later sections. Our first task is to
prove that $G_{A}\in \F_{A}$. In Section~2 we show  how this
follows from the facts that $\tilde G_\alpha\in \tilde \F_\alpha$
for all $\alpha:X\to [0,+\infty)$, $G_A=\tilde G_{\tilde\nu_A}$ if
$A$ is an effective divisor, and a variant of the Hironaka
desingularization theorem. By the same desingularization technique
we establish a representation of the Green function as the lower
envelope of the analytic disc functional $f\mapsto G_{f^*A}(0)$.
In Section~3 we study decomposition in ideal sheaves as a
preparation for Section~\ref{sec:properties} where we prove that
the estimates in the definition of the class $\F_A$ are locally
uniform.  This gives a direct proof of the relation $G_{A}\in
\F_{A}$ (without referring to desingularization), which in turn
implies certain refined maximality properties of the Green
function. In Section~5 we get a representation for the current
$(dd^cG_A)^p$ in the case when the ideal sheaf $\I_A$ has global
generators, and in Section~6 we study the case when the space is
reduced. In Section~7 we prove the product property of Green
functions, and finally in Section~8 we give a few explicit
examples.

\section{Definitions and main results}\label{sec:def}

We shall always let $X$ be a complex manifold and assume that $X$
is connected.  We denote  by $\PSH(X)$ the class of all
plurisubharmonic functions on $X$ and by $\PSH^-(X)$ its subclass
of all non-positive functions. We take $-\infty\in \PSH(X)$ and
set $\nu_{-\infty}=+\infty$.  We let $\O_X$ denote the sheaf of
germs of locally defined holomorphic functions on $X$.  We let $A$
be a closed complex subspace of $X$, $\I_A=(\I_{A,x})_{x\in X}$ be
the associated coherent sheaf of ideals in  $\O_X$, and $|A|$ be
the analytic variety in $X$ defined as the common set of zeros of
the locally defined functions on  $X$ with germs in $\I_A$. If $U$
is an open subset of $X$, then we let $\I_{A,U}$ denote the space
of all holomorphic functions on $U$ with germs in $\I_A$. We let
$\D$ denote the open unit disc in the complex plane $\C$ and $\T$
denote the unit circle.  We let $\O(Y,X)$ denote the set of all
holomorphic maps from a complex manifold $Y$ into  $X$. A map in
$\O(\D,X)$ is called an {\it analytic disc}, and if it can be
extended to a holomorphic map in some neighbourhood of the closed
disc $\overline \D$ then it is said to be {\it closed}. The
collection of all closed analytic discs is denoted by
$\O(\overline \D,X)$.

\begin{definition} Given a complex subspace $A$ of a connected complex
manifold $X$, the class $\F_A$ consists of all functions
$u\in\PSH^-(X)$ such that for every point $a\in X$ there exist
local generators $\psi_1,\dots,\psi_m$ for $\I_A$ near $a$ and a
constant $C$ depending on $u$ and the generators with $u\leq
\log|\psi|+C$ near $a$.
\end{definition}

Observe that $-\infty\in \F_{A}$ for every $A$.

\begin{definition} The {\it pluricomplex Green function $G_A$
with singularities along} $A$ is the upper envelope of all the
functions from the class $\F_A$, i.e.,
$$
G_A(x) = \sup\set{u(x);\: u\in\F_A}, \qquad x\in X.
$$
\end{definition}


The  local estimate $u\leq \log|\psi|+C$ is independent of the
choice of generators, i.e., if we have another set of generators
$\psi'=(\psi_1',\dots,\psi_k')$, then $u\leq \log|\psi'|+C'$ for
some constant $C'$. Furthermore, in the definition of the class
$\F_A$, $\psi=(\psi_1,\dots,\psi_m)$ can be replaced by any
holomorphic $\xi=(\xi_1,\dots,\xi_l)$, defined near $a$ and
satisfying
$$
\log|\xi|+c_1\le \log|\psi|\le \log|\xi|+c_2,
$$
which means precisely that the integral closure of the ideal
generated by the germs of the functions $\xi_i$ at $x$ coincides
with the integral closure of the ideal $\I_{A,x}$ for all $x$ in
some neighbourhood of $a$. (See \cite{D3}, Ch.~VIII, Cor.~10.5.)
We occasionally write $\log|\xi| \asymp\log|\psi|$ when
inequalities of this kind hold.

Let $X$ and $Y$ be complex manifolds and $\Phi:Y\to X$ be a
holomorphic map.  If $A$ is a complex subspace of $X$, then we
have a natural definition of a pullback $\Phi^*A$ of $A$ as a
complex subspace of $Y$.  The ideal sheaf $\I_{\Phi^*A}$ is
locally generated at a point $b$ by
$\Phi^*\psi_1,\dots,\Phi^*\psi_m$ if $\psi_1,\dots,\psi_m$ are
local generators for $\I_A$ at $\Phi(b)$. It is evident that
$\Phi^*u\in \F_{\Phi^*A}$ for all $u\in \F_A$, so
\begin{equation}
\Phi^*G_A\leq G_{\Phi^*A}. \label{eq:2.1}
\end{equation}

If $\Phi$ is proper and surjective and $v:Y\to \R\cup
\set{-\infty}$ is an upper semi-continuous function, then the
push-forward $\Phi_*v$ of $v$ to $X$ is  well defined by the
formula
$$
\Phi_*v(x)=\max_{y\in \Phi^{-1}(x)}v(y), \qquad x\in X.
$$

\begin{proposition}\label{propermap} Let $X$ and $Y$ be
complex manifolds of the same dimension and  $\Phi:Y\to X$ be a
proper surjective holomorphic map (for example, a finite branched
covering). Then $\Phi_*v\in\PSH(X)$ for all  $v\in\PSH(Y)$.
\end{proposition}

\begin{proof} In order to show that $\Phi_*v$ is upper
semicontinuous, we need to prove the relation $\Phi_*v(a)\geq
\limsup_{x\to a}\Phi_*v(x)$ for every $a\in X$.  We take a
sequence $a_j\to a$ such that $\Phi_*v(a_j)\to \limsup_{x\to
a}\Phi_*v(x)$.  Since $v$ is upper semicontinuous and $\Phi$ is
proper, there exist $b_j\in \Phi^{-1}(a_j)$ such that
$v(b_j)=\Phi_*v(a_j)$. By replacing $(b_j)$ by a subsequence we
may assume that $b_j\to b\in Y$. Then $\Phi(b)=a$ and
$$
\Phi_*v(a)\geq v(b)\geq \limsup_{j\to +\infty}v(b_j)=\lim_{j\to
+\infty} \Phi_*v(a_j)=\limsup_{x\to a}\Phi_*v(x).
$$
We let $V$ denote the set of all points $y$ in $Y$ for which
$d_y\Phi$ is degenerate.  Then $V$ is an analytic variety in $Y$
and Remmert's proper mapping theorem implies that $W=\Phi(V)$ is
an analytic variety in $X$.  It is sufficient to show that
$\Phi_*v$ is plurisubharmonic in a neighbourhood of every point
$a\in X\setminus W$, for the upper semicontinuity of $\Phi_*v$
then implies that $\Phi_*v\in \PSH(X)$.

Since $\Phi$ is a local biholomorphism on $Y\setminus
\Phi^{-1}(W)$, it follows that the fiber  $\Phi^{-1}(a)$ is
discrete and compact, thus finite, say that it consists of the
points $b_1,\dots,b_m$.  We choose a neighbourhood $U$ of $a$ in
$X\setminus W$ and biholomorphic maps $F_j:U\to F_j(U)\subseteq
Y\setminus \Phi^{-1}(W)$ with $F_j(a)=b_j$.  Then
$\Phi_*v(x)=\sup_{1\leq j\leq m}v\circ F_j(x)$ for all $x\in U$,
which shows that $\Phi_*v$ is plurisubharmonic in $U$.
\end{proof}

\bigskip
It is obvious that $u=\Phi_*\Phi^*u$ for all $u\in \PSH(X)$ and
$v\le \Phi^*\Phi_*v$ for all $v\in \PSH(Y)$.

\begin{proposition}\label{proper}
Let $X$ and $Y$ be complex manifolds of the same dimension, $A$ be
a closed complex subspace of $X$, and  $\Phi:Y\to X$ be a proper
surjective holomorphic map. Then $\Phi_*v\in \F_A$ for all $v\in
\F_{\Phi^*A}$ and
$$
\Phi^*G_A = G_{\Phi^*A}.
$$
\end{proposition}

\begin{proof}
If $a\in X$ and $\psi_1,\dots,\psi_m$ are local generators for
$\I_A$ near $a$, then  $v\leq \Phi^*\log|\psi|+C$ in some
neighbourhood of the compact set $\Phi^{-1}(a)$, which implies
$\Phi_*v\leq \log|\psi|+C$ near $a$. Hence we conclude from
Prop.~\ref{propermap} that $\Phi_*v\in \F_A$. Since $\Phi^*G_A\leq
G_{\Phi^*A}$, it is sufficient to prove that $v\leq \Phi^*G_A$ for
every $v\in \F_{\Phi^*A}$.  We have $\Phi_*v\in \F_A$, so $v\leq
\Phi^*\Phi_*v\leq \Phi^*G_A$.
\end{proof}

\bigskip
Our first main result is

\begin{theorem}\label{GinF}  If $X$ is a complex manifold and
$A$ is a closed complex subspace of $X$, then $G_A\in \F_A$.
\end{theorem}

Observe that in our definition of the class $\F_A$, the constant
$C$ in the local estimates $u\leq \log|\psi|+C$ is allowed to
depend both on the function $u$  and the local generators.  The
main work in our proof of Theorem~\ref{GinF} in
Sections~\ref{sec:decomposition} and \ref{sec:properties} is  to
prove that these estimates are indeed locally uniform, i.e., we
show that if $U$ is the domain of definition of $\psi$ and $K$ is
a compact subset of $U$, then there exists a constant $C_K$, only
depending on $K$ and $\psi$,  such that $u\leq \log|\psi|+C_K$ on
$K$ for all $u\in \F_A$. (See Lemma \ref{lub-lemma}.)

Let us now show how Theorem~\ref{GinF} follows from the facts that
$\tilde G_\alpha\in \tilde \F_\alpha$ for all $\alpha:X\to
 [0,+\infty)$, $G_A=\tilde G_{\tilde\nu_A}$ if $A$ is an effective
divisor, and the following variant of the Hironaka
desingularization theorem. (See \cite{BM}, Theorems 1.10 and
13.4.)

\medskip\noindent
{\em Given a closed complex subspace $A$ on a manifold $X$, there
exists a complex manifold $\hat X$ and a proper surjective
holomorphic map $\Phi:\hat X\to X$ which is an isomorphism outside
$\Phi^{-1}(|A|)$ and such that $\hat A=\Phi^*A$ is a
normal-crossing principal ideal sheaf (i.e., generated locally by
a monomial in suitable coordinates).}

\medskip\noindent
If we let $\Phi$ denote the desingularization map, then
$$
G_A=\Phi_*\Phi^*G_A=\Phi_*G_{\Phi^*A}=\Phi_*G_{\hat A}
=\Phi_*\tilde G_{\nu_{\hat A}}.
$$
Since $\tilde G_{\tilde\nu_{\hat A}}\in\F_{\hat A}$,
Proposition~\ref{proper} gives $G_A\in \F_A$ and the theorem is
proved.

If $X$ is one-dimensional, i.e., a Riemann surface, then  $\I_{A}$
is a principal ideal sheaf. If $\I_A=0$, the zero sheaf, then
$\tilde\nu_A=+\infty$. If $\I_A\neq 0$, then $|A|$ is discrete,
for $a\not\in |A|$ we have $\I_{A,a}=\O_{X,a}$ and
$\tilde\nu_A(a)=0$, and for  $a\in |A|$ the ideal $\I_{A,a}$ is
generated by the germ of $\zeta_a^m$ at $a$, where
$m=\tilde\nu_A(a)>0$ and $\zeta_a$ is a local generator for $\I_A$
near $a$ with $\zeta_a(a)=0$.  We obviously have
$$
G_A\geq \sum_{a\in \D}\tilde\nu_A(a)G_X(\cdot,a) \in \F_A,
$$
where $G_X(\cdot,a)$ is the Green function on $X$ with single pole
at $a$. In the special case $X=\D$, every function  $u\in
\F_A\setminus\set{-\infty}$ can be represented by the
Poisson--Jensen formula
$$
u(z)=\dfrac 1{2\pi}\int_\D G_\D(z,\cdot)\, \Delta u +\int_\T
P_\D(z,t)\, d\lambda_u(t), \qquad z\in \D,
$$
where $P_\D$ is the Poisson kernel for the unit disc $\D$ and
$\lambda_u$ is a nonpositive measure on the unit circle $\T$ (the
boundary value of $u$). We have $\nu_u(a)=\Delta u(\set a)/2\pi$,
so $\Delta u \geq 2\pi\sum_{a\in \D} \tilde\nu_A(a)\delta_a$,
where $\delta_a$ is the Dirac measure at the point $a$. Thus the
Poisson--Jensen formula implies
$$
u(z)\leq \int_\D G_\D(z,\cdot)\,\bigg(\sum_{a\in \D}
\tilde\nu_A(a)\delta_a\bigg)= \sum_{a\in \D}\tilde\nu_A(a)
G_\D(z,a)
$$
and we conclude that
$$
G_A(z)=\sum_{a\in \D} \tilde\nu_A(a)G_\D(z,a)= \sum_{a\in \D}
\tilde\nu_A(a)\log\bigg|\dfrac{z-a}{1-\bar az}\bigg|, \qquad z\in
\D,
$$
for every closed complex subspace $A$ of $\D$.

Now we let $X$ be any manifold, $f\in \O(\D,X)$ be an analytic
disc, and $a\in \D$. If $\psi_1,\dots,\psi_m$ are local generators
for $A$ at $f(a)$, then  $f^*\psi_1,\dots,f^*\psi_m$ are local
generators for $\I_{f^*A}$ near $a$.  If all these functions are
zero in some neighbourhood of $A$, then $\I_{f^*A}=0$ and
$\tilde\nu_{f^*A}=+\infty$. If one of them is not zero at $a$,
then $\I_{f^*A,a}=\O_{X,a}$ and $\tilde\nu_{f^*A}(a)=0$, and if
they have a common isolated zero at $a$, then
$\tilde\nu_{f^*A}(a)$ is the smallest positive multiplicity of
them.   Since $f^*G_A\leq G_{f^*A}$, we get
$$
G_A(x)\leq G_{f^*A}(0)=\sum_{a\in \D}\tilde\nu_{f^*A}(a)\log|a|,
\qquad f\in \O(\D,X),\quad f(0)=x.
$$

\begin{theorem}\label{th:envelope}
Let $X$ be a complex manifold and $A$ be a closed complex subspace
of $X$.  Then
$$
G_A(x)=\inf\set{G_{f^*A}(0)\, ;\, f\in \O(\overline \D,X),
x=f(0)}, \qquad x\in X.
$$
\end{theorem}

Let us show how the theorem follows from Hironaka's
desingularization theorem. If we  use the disc formula
(\ref{eq:disc}) for $\tilde G_\alpha$ with $\alpha=\tilde\nu_{\hat
A}$, the fact that $G_{\hat A}=\tilde G_\alpha$, and  the
desingularization map $\Phi$ above with $x=\Phi( \hat x)$, then
\begin{align*}
G_A(x)=\Phi^*G_{A}(\hat x)=G_{\hat A}(\hat x) &=\inf\set{\tilde
G_{g^*\alpha}(0)\,;\, g\in\O(\overline \D,\hat X),
g(0)=\hat x} \\
&\geq \inf\set{G_{g^*\hat A}(0)\,;\, g\in\O(\overline \D,\hat X),
g(0)=\hat x} \\
&= \inf\set{G_{g^*\Phi^*A}(0)\,;\, g\in\O(\overline \D,\hat X),
g(0)=\hat x} \\
&= \inf\set{G_{(\Phi_*g)^*A}(0)\,;\, g\in\O(\overline \D,\hat X),
g(0)=\hat x} \\
&\geq \inf\set{G_{f^*A}(0)\,;\, f\in\O(\overline \D, X), f(0)=x}
 \end{align*}
and we have proved Theorem \ref{th:envelope}. We will prove this
theorem without reference to desingularization or the disc formula
for $\tilde G_\alpha$ in a separate paper.

Let $X_1$ and $X_2$ be complex manifolds, $A_1$ and $A_2$ be
closed complex subspaces of $X_1$ and $X_2$, respectively,
$X=X_1\times X_2$ be the product manifold of $X_1$ and $X_2$, and
$A=A_1\times A_2$ be the product space of $A_1$ and $A_2$. If
$a=(a_1,a_2)\in X$ and $\psi^1_1,\dots,\psi^1_k$ and
$\psi^2_1,\dots,\psi^2_l$ are local generators for $\I_{A_1}$ and
$\I_{A_2}$ near $a_1$ and $a_2$, respectively, then the functions
$$
x=(x_1,x_2)\mapsto
\psi^1_1(x_1),\dots,\psi^1_k(x_1),\psi^2_1(x_2),\dots,\psi^2_l(x_2).
$$
are generators for $\I_A$ near $a$.  This implies that $X\ni
x=(x_1,x_2)\mapsto \max\set{u_1(x_1), u_2(x_2)}$ is in $\F_A$ for
all $u_1\in \F_{A_1}$ and $u_2\in \F_{A_2}$ , so we obviously have
$$
G_A(x) \geq \max\set{G_{A_1}(x_1), G_{A_2}(x_2)}, \qquad
x=(x_1,x_2)\in X.
$$
The following is called the {\it product property} for Green
functions.

\begin{theorem}\label{th:product} Let $X_1$ and $X_2$ be complex
manifolds, $A_1$ and $A_2$ be closed complex subspaces of $X_1$
and $X_2$, respectively, and $A$ be the product of $A_1$ and $A_2$
in $X=X_1\times X_2$.  Then
$$
G_A(x)=\max\set{G_{A_1}(x_1),G_{A_2}(x_2)}, \qquad x=(x_1,x_2)\in
X.
$$
\end{theorem}

We base our proof on  Th.~\ref{th:envelope} and give it in Section
\ref{sec:product}.

It was shown in \cite{LarSig2}, Prop.~3.2, that if $A$ is given by
a single holomorphic function with  effective divisor $Z_A$, then
$G_A$ satisfies $dd^cG_A\ge Z_A$ and, moreover, it is the largest
negative plurisubharmonic function with this property. (See the
last statement of Th.~3.3 in \cite{LarSig2}). Here $d=\partial +
\bar\partial$, $d^c= ( \partial -\bar\partial)/2\pi i$.

In the general case, a space $A$ generates  holomorphic chains
\begin{equation}\label{eq:reschain}
Z_A^p=\sum_i m_{i,p}[A_i^p],\end{equation}
where $A_i^p$ are
$p$-codimensional components of $|A|$ and $m_{i,p}\in\Z^+$.
Namely, if on a domain $U\subset X$ the space $A$ is given by
functions $\psi_1,\dots,\psi_m$  and $\codim\,|A|=p$ there, then
by the King-Demailly formula (\cite{D1}, Th.~6.20),
$$
(dd^c\log|\psi|)^p=\sum_i m_{i,p}[A_i^p] +R\quad{\rm on}\ U,
$$
where $m_{i,p}$ is the generic multiplicity of $\psi$ along
$A_i^p$ and $R$ is a positive closed current of bidegree $(p,p)$
on $U$, such that $\chi_{|A|}R=0$ and $\codim\,E_c(R)>p$ for every
$c>0$. Here $\chi_S$ is the characteristic function of a set $S$,
$E_c(R)=\{x\,;\, \nu_R(x)\ge c\}$ and $\nu_R(x)$ is the Lelong
number of the current $R$ at $x$. In other words, the holomorphic
chain $Z_A^p$ given by (\ref{eq:reschain}) is the residual
Monge-Amp\`ere current of $\log|\psi|$ on $|A|\cap U$.

\begin{theorem}\label{th:current}
Let $A$ have bounded global generators $\psi$ in $X$. Then
\begin{description}
\item{(i)} $G_A=\log|\psi|+O(1)$ locally near $|A|$. \item{(ii)}
If $\codim\,|A|=p$ on $U\subset X$, then  $(dd^cG_A)^p=Z_A^p +Q$
on $U$, where $Q$ is a positive closed current of bidegree $(p,p)$
on $U$, such that $\chi_{|A|}Q=0$ and $\codim\,E_c(Q)>p$ for every
$c>0$. If $U\cap |A|\subset J^p$, then $Q$ has zero Lelong
numbers; here the set $J^p$ consists of all points $a\in |A|$ such
that $p$ is the minimal number of generators of a subideal of
$\I_{A,a}$ whose integral closure is equal to  the integral
closure of $\I_{A,a}$.
\end{description}
\end{theorem}

A proof is given in Section \ref{sec:bounded} (and the sets $J^p$
are introduced and studied in Section~\ref{sec:decomposition}).

\section{Decomposition in ideal sheaves}\label{sec:decomposition}

In the case of complete intersection, i.e., when for every $a\in
|A|$ the local ideal $\I_{A,a}$ is generated by precisely
$p=\codim_a|A|$ germs of holomorphic functions, the relation
$G_A\in\F_A$ is in fact quite easy to prove without using the
desingularization technique. The main result of this section,
Prop.~\ref{decomp-theo}, gives a tool for the reduction of the
general situation to the complete intersection case in
Section~\ref{sec:properties}. Our approach develops a method from
\cite{R}.

We recall some basics on complex Grassmannians. (See, e.g.,
\cite{Chirka}, A3.4-5.) The Grassmannian $G(k,m)$ is the set of
all $k$-dimensional linear subspaces of $\C^m$ with the following
complex structure.  Let $S_{1\dots k}$ be the set of all $L\in
G(k,m)$ whose projections to the coordinate plane $\C_{1\dots k}$
of the variables $z_1,\dots,z_k$ are bijective.   Choosing a basis
$\set{(e_j,w_j)}$ in $L\in S_{1\dots k}$ with $e_j$ the standard
basis vectors in $\C^k$ and $w_j$ vectors in $\C^{m-k}$, we get a
representation of $L$ as the $k\times m$-matrix $(E,W)$, where $E$
is the unit $k\times k$-matrix and $W$ is a $k\times
(m-k)$-matrix.  This gives a parametrization of $S_{1\dots k}$ by
$k\times (m-k)$-matrices $W$.  In a similar way we parametrize all
the charts $S_I$, $I=(i_1,\dots,i_k)$.  Since the neighbouring
relations are holomorphic, this determines a complex structure on
$G(k,m)$.  It is easy to see that $\dim\,  G(k,m)=k(m-k)$.  The
set $\set{(z,L)\,;\, z\in L}\subset \C^m\times G(k,m)$ is
sometimes called the {\it incidence manifold}.

Let $\psi: \Omega\to\C^m$, $m>1$, be a holomorphic map on a domain
$\Omega$ in $\C^n$ and $Z=\{x\in \Omega\,;\,\psi(x)=0\}$. If $U$
is a subdomain of $\Omega$, then the graph
$\Gamma_U=\{(x,\psi(x))\,;\,x\in U\}$ of $\psi$ over $U$ is an
$n$-dimensional complex manifold in $\C^n\times\C^m$. Given $k\le
m-1$, let $\Gamma_U^k$ be the pullback of $\Gamma_U$ to the
incidence variety in $\Gamma_U\times G(k,m)$. Namely, $\Gamma_U^k$
is the closure of the set
$$
\{(x,\psi(x),L)\,;\,  x\in U\setminus Z,\ L\in G(k,m),\ \psi(x)\in
L\}.
$$
By $\rho_k$ we denote the projection from $\Gamma^k_\Omega$ to
$G(k,m)$, and by $\pi_k$ its projection to $\Omega$.

For $x\in \Omega\setminus Z$, the fiber $\rho_k\circ\pi_k^{-1}(x)$
consists of all $L\in G(k,m)$ passing through $\psi(x)\neq 0$ and
thus is isomorphic to $G(k-1,m-1)$. Therefore $\dim\,
\Gamma_\Omega^k=n+(k-1)(m-k)$.

Let $I^k=I^k(\psi)$ be the collection of all points $x$ in
$\Omega$ such that $\rho_k(\Gamma_U^k)=G(k,m)$ for every
neighbourhood $U$ of $x$, i.e.,
$\rho_k\circ\pi_k^{-1}(x)=G(k,m)$. Evidently, $I^1\subseteq
I^2\subseteq\ldots\ \subseteq I^{m-1}\subseteq Z$.

\begin{lemma}\label{dimind-lemma} $I^k$ is an analytic set of
dimension at  most $n-m+k-1$.
\end{lemma}

\begin{proof} We have
\begin{equation}\label{eq:prod} \pi_k^{-1}(I^k)=I^k\times\{0\}\times G(k,m),
\end{equation}
so $I^k\subset \pi_k\circ\rho_k^{-1}(L)$ for each $L\in G(k,m)$.
On the other hand, for every $x\not\in I^k$ there exists $L\in
G(k,m)$ such that $x\not\in\pi_k\circ\rho_k^{-1}(L)$. Thus
$$I^k=\bigcap_{L\in G(k,m)} \pi_k\circ\rho_k^{-1}(L).
$$
Each $\rho_k^{-1}(L)$ is an analytic set in $\Gamma^k_\Omega$.
Since the map $\pi_k$ is proper, Remmert's theorem implies that
$\pi_k\circ\rho_k^{-1}(L)$ is an analytic subset of $\Omega$ for
any $L$, and so is $I^k$.

The set $\pi_k^{-1}(I^k)$ is a nowhere dense analytic subset of
$\Gamma^k_\Omega$, and thus
$\dim\,\pi_k^{-1}(I^k)<\dim\,\Gamma^k_\Omega=n+(k-1)(m-k)$. By
(\ref{eq:prod}), $\dim\,\pi_k^{-1}(I^k)=\dim\,I^k+k(m-k)$.
Therefore $\dim\,I^k<n+(k-1)(m-k)-k(m-k)=n-m+k$.
\end{proof}

\begin{corollary}\label{cor:empty}
If $m>n$, then $I^k=\emptyset$ for all $k\le m-n$.
\end{corollary}

\begin{lemma}\label{reduce-lemma}
For any $a\in Z\setminus I^k$ there exist a neighbourhood $U$ of
$a$ and holomorphic functions $\xi_1,\ldots,\xi_{m-k}$ (linear
combinations of $\psi_1,\ldots,\psi_m$) such that
$\log|\psi|\asymp \log|\xi|$ in $U$.
\end{lemma}

\begin{proof}
Given $a\in Z\setminus I^k$, one can find a neighbourhood $U$ of
$a$ such that $\rho_k(\Gamma_U^k)\neq G(k,m)$. Since the set
$G(k,m)\setminus \rho_k(\Gamma_U^k)$ is open, there exists $L_0$
in the chart $S_{1\ldots k}$ of $G(k,m)$ such that
\begin{equation}\label{eq:avoid}
\psi(x)\cap\omega=\emptyset
\end{equation}
 for some neighbourhood $\omega\subset
S_{1\ldots k}$ of $L_0$ and all $x\in U\setminus Z$.

Let $(E,W_0)$ be the canonical representation of $L_0$. For every
$y=(y',y'')\in \C^k\times\C^{m-k}$ with $y'\neq 0$, the map
$y\mapsto (y',y'W_0)$ is the projection to the space $L_0$. By
elementary linear algebra arguments (see Lemma~\ref{linear-lemma}
below), relation (\ref{eq:avoid}) implies existence of $r>0$ such
that
\begin{equation}\label{eq:elbound}
|\psi''(x)-\psi'(x)W_0|\ge r|\psi'(x)|, \quad x\in U.
\end{equation}
We define a map $\xi:U\to\C^{m-k}$ by
$\xi(x)=\psi''(x)-\psi'(x)W_0$. Then
$$
|\xi(x)|\le C|\psi(x)|,\quad x\in U.
$$
Furthermore, inequality (\ref{eq:elbound}) implies
$$
|\psi|^2\le |\psi'|^2+2|\psi''-\psi'W_0|^2+2|\psi'W_0|^2\le
C|\psi''-\psi'W_0|^2=C|\xi|^2,
$$
and the assertion follows.
\end{proof}

\begin{lemma}\label{linear-lemma} Let $W_0$ be a complex $k\times
(m-k)$-matrix and a set $S\subset\C^k\times\C^{m-k}$ be such that
$|y''-y'W|>0$ for all $y=(y',y'')\in S$ and all matrices $W\in
\C^{k(m-k)}$ with $|W-W_0|<\delta$ (all the norms $|\cdot|$ are
the Euclidean norms in the corresponding linear spaces). Then
$$|y''-y'W|\ge\frac\delta{k}|y'|,\quad y\in S,\ |W-W_0|<\delta.$$
\end{lemma}

\begin{proof} Suppose there exists $y\in S$ and $W$ in the
$\delta$-neighbourhood of $W_0$ such that
$|y''-y'W|<\frac\delta{k}|y'|$. For the vector $z=(z',z''):=
(y',y''-y'W_0)$ this means $|z''|<\tfrac{\delta} {k}|z'|$.

We choose $l\in [1,k]$ such that $|z_l|=\max\{|z_i|\,;\, 1\le i\le
k\}$ and consider the $k\times (m-k)$-matrix $V$ with the entries
$V_{lj}=z_{k+j}/z_l$ for $1\le j\le m-k$, and $V_{ij}=0$ for all
$i\neq l$ and  $1\le j\le m-k$. Then
$$
|V|=\frac{|z''|}{|z_l|}\le \frac{|z''|}{k|z'|}<\delta
$$
and $z'V=z''$. The latter relation is equivalent to $y''-y'W=0$
with $W=W_0+V$. Since $|W-W_0|=|V|<\delta$, this contradicts the
hypothesis of the lemma.
\end{proof}

\bigskip
We recall that the {\it analytic spread} of an ideal $\I$ equals
the minimal number of generators of a  subideal of $\I$ whose
integral closure coincides with the integral closure of $\I$, see
\cite{NR}.

\begin{proposition}\label{decomp-theo} Let $A$ be a closed complex subspace of
a manifold $X$, $\dim\,X=n$. Then the set $|A|$ can be decomposed
into the disjoint union of  local analytic varieties $J^k$, $1\le
k\le n$, such that
\begin{description} \item{(i)}
$\codim\,J^k\ge k$ and \item{(ii)} for each $a\in J^k$, the ideal
$\I_{A,a}$ has analytic spread at most $k$.
\end{description}
\end{proposition}

\medskip
\begin{proof}
Let $\psi=(\psi_1,\ldots,\psi_m)$ be generators of $\I_A$ on a
domain $\Omega\subset X$. Set $N=\min\,\{n,m\}$,
$Z=|A|\cap\Omega$, $J^1=Z\setminus I^{m-1}$,
$J^k=I^{m-k+1}\setminus I^{m-k}$ for $k=2,\ldots,N-1$, and
$J^{N}=I^{m-N+1}$ (some of them can be empty). The sets $J^k$ are
pairwise disjoint, $\dim\,J^k\le n-k$ (Lemma~\ref{dimind-lemma}),
and $\cup_k J^k=Z$. On a neighbourhood of each point of $J^k$, the
singularity of the function $\log|\psi|$ is equivalent to one
defined by the function $\log|\xi|$ with
$\xi=(\xi_1,\ldots,\xi_k)$ (this follows from
Lemma~\ref{reduce-lemma}, if $m\leq n$, and
Corollary~\ref{cor:empty}, in the case $m>n$). This means that the
ideal generated by the germs of $\psi_i$ at $a\in J^k$ has
analytic spread at most $k$.

Let $\psi'=(\psi_1',\ldots,\psi_{m'}')$ be other generators of
$\I_A$ on $\Omega$; by adding some identically zero components to
either $\psi$ or $\psi'$ we can assume $m'=m$. For any point $a\in
Z\setminus I^k(\psi)$, relation (\ref{eq:avoid}) implies existence
of a neighbourhood $U'$ of $a$ and a plane $L_0'\in G(k,m)$ such
that $\psi'(x)\cap\omega'=\emptyset$ for some neighbourhood
$\omega'$ of $L_0'$ and all $x\in U'\setminus Z$, so $a\in
Z\setminus I^k(\psi')$. This shows that the sets $J^k$ are
independent of the choice of generators of $\I_{A,\Omega}$.
Therefore each $J^k$ is well defined as a local (not necessarily
closed) analytic variety in $X$ with properties (i) and (ii).
\end{proof}

\begin{example}
Let $A$ be generated by $\psi(x)=(x_1^2x_2,x_1^2x_3,x_1x_2x_3)$ in
$\C^3$. Then $|A|=\C_{23}\cup \C_1$; here $\C_{23}$ is the
coordinate plane of the variables $x_2$ and $x_3$, i.e.,
$\C_{23}=\set{x_1=0}$, and $\C_1=\set{x_2=x_3=0}$. The variety
$|A|$ has the decomposition $|A|=J^1\cup J^2\cup J^3$ with
$J^1=\C_{23}\setminus (\C_2\cup\C_3)$, $J^2=\C_1^*\cup\C_2^*\cup
\C_3^*$, and $J^3=\set{0}$. Near points of $J^1$ we have
$\log|\psi|\asymp \log|x_1|$. As to $J^2$, the relation
$\log|\psi|\asymp \log|\xi|$ is satisfied with $\xi=(x_2,x_3)$
near points of $\C_1^*$, and we can take $\xi=(x_1^2,x_1x_3)$ near
points of $\C_2^*$ and $\xi=(x_1^2,x_1x_2)$ near points of
$\C_3^*$.
\end{example}

\section{Upper bounds and maximality}\label{sec:properties}

We recall that a function $u\in \PSH(X)$ is called {\it maximal}
in $X$ if for every relatively compact subset $U$ of $X$ and for
each upper semicontinuous function $v$ on $\overline U$ such that
$v\in \PSH(U)$ and $v\leq u$ on $\partial U$, we have $v\leq u$ in
$U$. An equivalent form is that for any $v\in \PSH(X)$ the
relation $\set{v>u}\Subset X$ implies $v\leq u$ on $X$.

We will use the following variant of the maximum principle for
unbounded plurisubharmonic functions.

\begin{lemma}\label{maxprin-lemma}
Let $D\subset\C^k$ be a bounded domain and  $u, v\in\PSH(D)$ such
that

\begin{description}
\item{(i)} $v$ is bounded above,

\item{(ii)} the set $S:=v^{-1}(-\infty)$ is closed in $D$,

\item{(iii)} $v$ is locally bounded and maximal on $D\setminus S$,

\item{(iv)} for any $\epsilon>0$ there exists a compact
$K_\epsilon \subset D$ such that $u(z)\le v(z)+\epsilon$ on
$D\setminus K_\epsilon$, and

\item{(v)} $\displaystyle \limsup_{z\to a,\:z\not\in
S}(u(z)-v(z))<\infty$ for each $a\in S$.
\end{description}
\noindent Then $u\le v$ in $D$.
\end{lemma}

\begin{proof} By (i) we may assume that $v$ is negative in $D$.
Take any $\epsilon>0$ and $\delta>0$.  Then it is sufficient to
prove that $u_1=(1+\delta)(u-\epsilon)\leq v$. By (v) we conclude
that each point $a\in S$ has a neighbourhood $U_a\Subset D$ where
$u_1\le v$ and by (iv) that there is a domain $D_1\Subset D$ such
that $u_1\le v$ on $D\setminus D_1$. By (ii) $S\cap \overline D_1$
is compact, so we can take a finite covering of $S\cap \overline
D_1$ by $U_{a_j}$, $1\leq j\leq N$. Then
$D_2=D_1\setminus\bigcup_j\overline U_{a_j}$ is an open subset of
$D$ on which  $v$ is bounded and $u_1\le v$ holds on $\partial
D_2$. By (iii) $v$ is maximal on $D_2$, so $u_1\leq v$ on $D_2$
and thus on $D$.
\end{proof}

\medskip
The next statement is the crucial point in the proof that
$G_A\in\F_A$.

\begin{lemma}\label{lub-lemma}
Let $\psi=(\psi_1,\ldots,\psi_m)$ be a holomorphic map on a domain
$\Omega\subset\C^n$ and $Z$ be its zero set. Then for every
$K\Subset \Omega$ there exists a number $C_K$ such that any
function $u\in PSH^-(\Omega)$ which satisfies $u\le \log|\psi| +
O(1)$ locally near points of $Z$ has the bound $u(x)\le
\log|\psi(x)| + C_K$ for all $x\in K$.
\end{lemma}

\begin{proof} What we need to prove is that each point $a\in Z$
has a neighbourhood $U$ where $u\le\log|\psi|+C$ with $C$
independent of the function $u$.

Let $\codim_a\,Z=p$. Then, by Prop.~\ref{decomp-theo}, $a\in
J^{k}$ for some $k\in [p,n]$ and thus there exist $k$ holomorphic
functions $\xi_1,\ldots,\xi_k$ such that $\log|\xi|\asymp
\log|\psi|$ near $a$. We will argue by induction in $k$ from $p$
to $n$.

Let $a\in J^{p}$; this means that there is a neighbourhood $V$ of
$a$ such that $Z\cap V$ is a complete intersection given by the
functions $\xi_1,\ldots,\xi_{p}$. By Thie's theorem \cite{Thie},
(see also \cite{D1}, Th.~5.8), there exist local coordinates
$x=(x',x'')$, $x'=(x_1,\ldots,x_p)$, $x''=(x_{p+1},\ldots,x_n)$,
centered at $a$ and balls $\B'\subset {\C}^p$, $\B''\subset
{\C}^{n-p}$ such that $\B'\times \B''\Subset V$, $Z\cap(\B'\times
\B'')$ is contained in the cone $\{|x'|\le \gamma |x''|\}$ with
some constant $\gamma>0$, and the projection of $Z\cap(\B'\times
\B'')$ onto $\B''$ is a ramified covering with a finite number of
sheets. Let $r_1=2\gamma r_2$ with a sufficiently small $r_2>0$ so
that $\B_{r_1}'\subset \B'$ and $\B_{r_2}''\subset \B''$, then for
some $\delta>0$
$$|\xi(x)|\ge\delta,\qquad  x\in \partial \B_{r_1}'\times
\B_{r_2}''.
$$

Given $x_0''\subset \B_{r_2}''$, denote by $Z(x_0'')$ and
$\Sing\,Z(x_0'')$ the intersections of the set
$\B_{r_1}'\times\{x_0''\}$ with the varieties $Z$ and $\Sing\,Z$,
respectively. Since the projection is a ramified covering,
$Z(x_0'')$ is finite for any $x_0''\in \B_{r_2}''$, while
$\Sing\,Z(x_0'')$ is empty for almost all $x_0''\in \B_{r_2}''$
because $\dim\, \Sing\,Z\le n-p-1$; we denote the set of all such
generic $x_0''$ by $E$.

Fix any $x_0''\in E$ and consider the function
$$v(x')=\log(|\xi(x',x_0'')|/\delta).$$
It is plurisubharmonic on $\B_{r_1}'$, nonnegative on $\partial
\B_{r_1}'$ and maximal on $\B_{r_1}'\setminus Z(x_0'')$, since the
map  $\xi(\cdot,x_0''):\B_{r_1}'\to{\C}^p$ has no zeros outside
$Z(x_0'')$.

For any function $u\in PSH^-(Y)$ which satisfies $u\le \log|\xi| +
O(1)$ locally near regular points of $Z$, we have, by
Lemma~\ref{maxprin-lemma}, $u(x',x_0'')<v(x')$ on the whole ball
$\B_{r_1}'$.

Since $x_0''\in E$ is arbitrary, this gives us $u\le
\log|\xi|-\log\delta$ on $\B_{r_1}'\times E$. The continuity of
the function $\log|\xi|$ extends this relation to the whole set
$U=\B_{r_1}'\times \B_{r_2}''$, which proves the claim for $k=p$.

Now we make a step from $k-1$ to $k$. Since $\dim\, J^{k}\le n-k$,
we use Thie's theorem to get a coordinate system centered at $a\in
J^{k}$ such that the projection of $J^{k}\cap (\B'\times\B'')$ to
$\B''\subset\C^{n-k}$ is a finite map and
$(\partial\B'\times\overline{\B''})\cap J^i=\emptyset$ for all
$i\ge k$. Therefore, by the induction assumption and a compactness
argument, $u\le\log|\xi|+C$ near
$\partial\B'\times\overline{\B''}$, where the constant $C$ is
independent of $u$.

Now for any $x_0''\subset \B''$ we consider the function
$v(x')=\log|\xi(x',x_0'')|+C$. Then Lemma~\ref{maxprin-lemma}
gives us $u(x',x_0'')<v(x')$ on $\B'$ and hence $u\le\log|\xi|+C$
on $\B'\times\B''$.
\end{proof}

\medskip

{\it Remark.} Note that the uniform bound $u\le\log|\psi|+C$ near
points $a\in J^p$, $\codim_a Z =p$, was deduced from the local
bounds only near regular points of $Z$.

\medskip

\begin{prooftx}{Proof of Theorem~\ref{GinF}}
The relation $G_A\le\log|\psi|+O(1)$ follows from
Lemma~\ref{lub-lemma}. This implies that its upper semicontinuous
regularization $G_A^*$ is in $\F_A$ and thus $G_A^*=G_A$.
\end{prooftx}

One of the most important properties of the ``standard''
pluricomplex Green function $G_{X,a}$ with logarithmic pole at
$a\in X$ is that it satisfies the homogeneous Monge-Amp\`ere
equation $(dd^cG_{X,a})^n=0$ outside the point $a$; in other
words, $G_{X,a}$ is a maximal plurisubharmonic function on
$X\setminus\{a\}$. In our situation, one can say more.

\begin{theorem}\label{max-theo} The function  $G_A$ is maximal on
$X\setminus |A|$ and locally maximal outside a discrete subset of
$|A|$ (actually, the set $J^n$ from Prop.~\ref{decomp-theo}). If
$A$ has $k<n$ global generators on $X$, then $G_A$ is maximal on
the whole $X$.
\end{theorem}

\begin{proof}
Take any point $a\not\in J^n$. By Proposition~\ref{decomp-theo},
there exist functions $\xi_1,\ldots,\xi_{k}\in\I_{A,U}$, $k<n$,
generating an ideal whose integral closure coincides with the
integral closure of $\I_{A,U}$, and so $G_A\le\log|\xi|+C$ on $U$.
The function $\log|\xi|$ is maximal on $U$, which follows from the
fact that it is the limit of the decreasing sequence of maximal
plurisubharmonic functions $u_j=\frac 12\log(|\xi|^2+{\frac 1j})$.
(See \cite{R}, Example~1.) Take any domain $W\Subset U$. Given a
function $v\in PSH(U)$ with $v\le G_A$ on $U\setminus W$, we have
to show that $v\le G_A$ on $U$. Consider the function $w$ such
that $w=G_A$ on $X\setminus W$ and $w=\max\{G_A,v\}$ on $W$. Since
$G_A\le\log|\xi|+C$ on $U$, we have $w\le \log|\xi|+C$ on
$U\setminus W$, and the maximality of $\log|\xi|$ on $U$ extends
this inequality to the domain $W$. Therefore, $w\in\F_A$ and thus
$w\le G_A$ on $U$.

When $a\not\in |A|$, we can take $U=X\setminus |A|$ and $\xi\equiv
1$, which gives us maximality of $G_A$ on $U=X\setminus |A|$.

Finally, if $A$ has $k<n$ global generators on $X$, then the same
arguments with $U=X$ show the maximality of $G_A$ on the whole
$X$.
\end{proof}

\medskip

{\it Remark.} If $J^n=\emptyset$, the Green function is locally
maximal on the whole $X$. We don't know if this implies its
maximality on $X$.

\section{Complex spaces with bounded global generators}\label{sec:bounded}

If $A$ has bounded generators $\psi$, which we can choose such
that $|\psi|<1$, then $\log|\psi|\in\F_A$. This gives immediately

\begin{proposition}\label{th:asymp} Let $A$ be a closed complex subspace of a manifold
$X$ and assume that $A$ has bounded global generators $\{\psi_i\}$
(for example, $X$ is a relatively compact domain in a Stein
manifold $Y$ and $A$ is a restriction to $X$ of a complex space
$B$ on $Y$), then
\begin{equation}\label{eq:asymp}
 G_A=\log|\psi|+O(1)
\end{equation}
locally near $|A|$.
\end{proposition}

To describe the boundary behaviour of $G_A$, we recall the notion
of strong plurisubharmonic barrier. Let $X$ be a domain in a
complex manifold $Y$, and let $p\in \partial X$. A
plurisubharmonic function $v$ on $X$ is called a {\it strong
plurisubharmonic barrier at} $p$ if $v(x)\to 0$ as $x\to p$, while
$\sup_{X\setminus V}v <0$ for every neighbourhood $V$ of $p$ in
$Y$. By standard arguments (see, e.g., \cite{LarSig2},
Proposition~2.4) we get

\begin{proposition}
Let $X$ be a domain in a complex manifold $Y$, and let a closed
complex subspace of $X$ have bounded global generators. If $X$ has
a strong plurisubharmonic barrier at $p\in\partial
 X\setminus |A|$, then
$G_A(x)\to 0$ as $x\to p$.
\end{proposition}


A uniqueness theorem for the Green function is similar to that for
the divisor case in \cite{LarSig2}, but the proof is different
(since the function $u-\log|\psi|$ need not be plurisubharmonic)
and follows from Lemma~\ref{maxprin-lemma} and
Proposition~\ref{th:asymp}.

\begin{theorem}  Let a complex space $A$ have bounded global generators
$\psi_i$ on $X$, and let a function $u\in PSH^-(X)$ have the
properties
\begin{description}
\item{(i)} $u$ is locally bounded and maximal on $X\setminus |A|$.
\item{(ii)} For any $\epsilon>0$ there exists a compact subset $K$
of $X$ such that $u\ge G_A-\epsilon$ on $X\setminus K$;
\item{(iii)} $u=\log|\psi|+O(1)$ locally near $|A|$.
\end{description}
Then $u=G_A$.
\end{theorem}

Relation (\ref{eq:asymp}) allows us to derive the properties of
the Monge-Amp\'ere current $(dd^cG_A)^p$.

\begin{prooftx}{Proof of Theorem \ref{th:current}}
Since $G_A$ is locally bounded on $U\setminus |A|$ and
$\codim\,|A|=p$, the current $(dd^cG_A)^p$ is well defined on $U$.
Moreover, Siu's structural formula  for positive closed currents
\cite{Siu} (see also \cite{D1}, Theorem~6.19) gives us a (unique)
representation for the current $(dd^c G_A)^p$ as
$$(dd^c G_A)^p=\sum_j\lambda_j[B_j]+Q,
$$
where $B_j$ are some irreducible analytic varieties of codimension
$p$, $\lambda_j$ are the generic Lelong numbers of $(dd^c G_A)^p$
along $B_j$, i.e.,
$$
\lambda_j=\inf\{\nu((dd^c G_A)^p,a):a\in B_j\},
$$
and $Q$ is a positive closed current such that $
\codim\{x:\nu(Q,x)\ge c\}>p$ for each $c>0$.

As $G_A$ has asymptotics (\ref{eq:asymp}) near points of the set
$|A|$, Demailly's Comparison Theorem for Lelong numbers
(\cite{D1}, Theorem~5.9) implies
$$\nu((dd^c G_A)^p,a)=\nu((dd^c \log|\psi|)^p,a)$$
at every point $a\in |A|\cap J^p\cap U$. In particular, the
generic Lelong number of $(dd^c G_A)^p$ along each variety $A_i^p$
equals the multiplicity of this component in $|A|$. Besides,
$\nu((dd^c G_A)^p,a)=0$ for any $a\not\in |A|$. This shows that
$\{B_j\}_j$ are exactly the $p$-codimensional components of the
variety $|A|$ in $U$ and $\sum\lambda_j[B_j]=Z_A^p$ on $U$.

Finally, if $U\cap |A|\subset J^p$, then $U\cap|A|$ can be given
locally by $p$ holomorphic functions $\xi_i$ with
$\log|\xi|\asymp\log|\psi|$. By King's formula,
$(dd^c\log|\xi|)^p=Z_A^p$, which means, in particular, that
$(dd^c\log|\xi|)^p$ has zero Lelong numbers outside $\cup_i
A_i^p$. Since the currents $(dd^c G_A)^p$ and $(dd^c\log|\xi|)^p$
have the same Lelong numbers, this proves the last statement.
\end{prooftx}

\medskip
So the Green function satisfies, as in the divisor case, the
relation $(dd^cG_A)^p\ge Z_A^p$, but for $p>1$ it is not the
largest negative plurisubharmonic function with this property
(even for reduced spaces that are complete intersections).
 For example, let $X$ be
the unit polydisc in ${\C}^3$ and $A$ be generated by
$\psi(z)=(z_1,z_2)$. Then $G_A=\max\{\log|z_1|,\log|z_2|\}$ and,
moreover, $(dd^cG_A)^2=Z_A=[A]$.  But the functions
$u_N=\max\{N\log|z_1|,N^{-1}\log|z_2|\}$, $N>0$, also satisfy
$(dd^c u_N )^2=[A]$, although they are not dominated by $G_A$. It
is easy to see that the upper envelope of all such functions
equals $0$ outside $|A|$ and $-\infty$ on $|A|$. Therefore, in the
case $\codim\,|A|>1$ there is no counterpart for the description
of the Green function in terms of the current $Z_A$.

\section{Reduced spaces}\label{sec:reduced}

Now we return to relations between the functions $G_A$ and $\tilde
G_{\tilde\nu_A}$ (see Introduction). As was already mentioned, one
has always $G_A\le\tilde G_{\tilde\nu_A}$ and $G_A<\tilde
G_{\tilde\nu_A}$ for 'generic' spaces $A$, however $G_A=\tilde
G_{\tilde\nu_A}$ for effective divisors $A$. Here we show that the
equality holds also in the case of {\sl reduced} complex spaces.

When $A$ is a reduced space, it can be identified with the
analytic variety $|A|$. Its generators $\psi_1,\ldots,\psi_m$ on
$U$ have the property: if a holomorphic function $\varphi$
vanishes on $A\cap U$, then $\varphi=\sum h_i\psi_i$ with
$h_i\in\O(U)$.

Since $\tilde\nu_A=1$ at all regular points of $A$, it is natural
to consider the class
$$
\tilde\F_A^1=\{u\in PSH^-(X);\: \nu_u(a)\ge 1 \text{ for all }
a\in Reg\,A\}.
$$
Note that upper semicontinuity of the Lelong numbers implies
$\nu_u\ge 1$ on the whole $A$.

We evidently have  $\tilde\F_A^1 \subseteq
\tilde\F_{\tilde\nu_A}\subseteq \F_A$.

\begin{theorem} If $A$ is a reduced subspace of $X$, then
$\tilde\F_A^1 = \tilde\F_{\tilde\nu_A}= \F_A$ and consequently
$$
G_A(x)=\tilde G_{\tilde\nu_A}(x)=\sup\set{u(x)\,;\,\ u\in
\tilde\F_A^1}.
$$
\end{theorem}

\begin{proof} It suffices to show that for any function $u\in\tilde\F_A^1$
and every point $a\in A$ there is a neighbourhood $U$ of $a$ and a
constant $C$ such that
\begin{equation}\label{eq:regloc}
u(x)\le \log|\psi(x)|+C, \qquad  x\in U.
\end{equation}

We will use induction on the dimension of $X$. The case
$\dim\,X=1$ is evident. Assume it proved for all $X$ with
$\dim\,X<n$ and take any $u\in \tilde\F_A^1$. When $\dim_a A=0$,
relation (\ref{eq:regloc}) follows easily from the fact that
$\log|\psi(x)|=\log|\zeta(x)|+O(1)$ near $a\in A$, where $\zeta$
are local coordinates near $a$ with $\zeta(a)=0$. So we assume
$\dim_a A>0$. We first treat the case when $a$ is a regular point
of $A$, $\codim_a A=p<n$. Since the problem is local, we may then
assume that $X\subset\Cn$ and contains the unit polydisc ${\D}^n$,
$a=0$, and the restriction $A'$ of $A$ to ${\D}^n$ is given by
$\psi(x)=(x_1,\ldots,x_p)$. Then the restriction of $u$ to
${\D}^n$ is dominated by the Green function $\tilde
G_{\tilde\nu_{A'}}$. By the product property for this type of
Green function (\cite{LarSig2}, Theorem~2.5), $\tilde
G_{\tilde\nu_{A'}}(x)=\max\{\log|x_j|,\ 1\le j\le p\}$. This
implies (\ref{eq:regloc}) for $a\in Reg\,A$.

For $a\in Sing\,A$ we will argue similarly to the proof of
Lemma~\ref{lub-lemma}. There is a neighbourhood $V$ of $a$ such
that $V\cap Sing\,A\subset J^p\cup J^{p+1}\cup\ldots J^n$. The
proof for $a\in J^{k}$, $p\le k\le n$, is then by induction in
$k$.

For $a\in J^p\cap V$ relation (\ref{eq:regloc}) follows directly
from the remark after Lemma~\ref{lub-lemma}.

Assuming (\ref{eq:regloc}) proved for $a\in J^p\cup\ldots\cup
J^{k}$, we take $a\in J^{k+1}$. We choose coordinates
$x=(x',x'')\in\C^{k+1}\times\C^{n-k-1}$ such that $a=0$, the
projection of $J^{k+1}\cap \B$ to $\B''$ is a finite map and
$\partial\B'\times\B''\cap J^i=\emptyset$ for all $i\ge k+1$, so
the $k$-induction assumption gives
\begin{equation}\label{eq:Lbound}
u(x)\le \log|\psi(x)|+C,\quad x\in\partial\B'\times\B''.
\end{equation}
 Take any $b=(b',b'')\in
\B'\times\B''$ and consider the $(k+1)$-dimensional plane
$L=\{x\,;\,x''=b''\}$. Then the restriction $u_L$ of $u$ to the
plane $L$ (in the same way we will use the denotation $\psi_L$,
$\B_L$, $A_L$, etc.) has Lelong numbers at least $1$ at all points
of $A_L$, so $u_L\in\tilde\F_{A_L,\B_L}^1$.  Since $\dim\,\B_L<n$
and the components of $\psi_L$ generate $A_L$, the $n$-induction
assumption implies $u_L\in\F_{A_L,\B_L}$. Therefore, $u_L\le
\log|\psi_L|+O(1)$ locally near points of $A_L$.

Since $a\in J^{k+1}$, we can find functions
$\xi_1,\ldots,\xi_{k+1}$ such that $\log|\xi|\asymp\log|\psi|$ on
$\B$. Therefore $u_L\le \log|\xi_L|+O(1)$ locally near all points
of $A_L$, and, by (\ref{eq:Lbound}), $u_L\le \log|\psi_L|+C_1$ on
a neighbourhood of $\partial \B_L$ with $C_1$ independent of $L$.
The function $\xi_L$ is maximal on $\B_L\setminus A_L$, so by
Lemma~\ref{maxprin-lemma}, $u_L\le \log|\xi_L|+C_1$ everywhere on
$\B_L$. Since the plane $L$ was chosen arbitrary, this gives us
(\ref{eq:regloc}) for $a\in J^{k+1}$.

This proves the inductive step in the induction in $k$ and, at the
same time, in the induction in $n$.
\end{proof}

\medskip

\medskip

Theorem~\ref{max-theo} for reduced spaces has the following form
(compare with the remark after the proof of
Theorem~\ref{max-theo}).

\begin{theorem} The Green function of a reduced space
$A$ is maximal on $X\setminus A_0$, where $A_0$ is the collection
of $0$-dimensional components of $A$.
\end{theorem}

\begin{proof} We need to show that for every domain
$U\Subset X':=X\setminus A_0$ and a function $u\in PSH(X')$ the
condition $u\le G_A$ on $X'\setminus U$ implies $u\le G_A$ on $U$.

Consider the set $E_1(u)=\{x\in X\,;\, \nu_u(x)\ge 1\}$. Since
$u\le G_A$ on $X'\setminus U$, we have $E_1(u)\setminus U\supset
A\setminus U$. By Siu's theorem, $E_u$ is an analytic variety in
$X$, so it must contain the whole $A$. This means that
$u\in\tilde\F_A^1$ and thus is dominated by $\tilde
G_{\tilde\nu_A}=G_A$ on $X$.
\end{proof}

\section{The  product property}\label{sec:product}

Our proof of Th.~\ref{th:product} in this section is based on
Th.~\ref{th:envelope}. It is a modification of the proof of
Th.~2.5 in \cite{LarSig2} which in turn generalizes a proof of
Edigarian \cite{E} of the product property for the single pole
Green function. For the sake of completeness we have repeated some
arguments from \cite{E}, \cite{LarSig2}, and \cite{LarLasSig1}.

We introduce the following notation:  If the function $\varphi$ is
holomorphic in some neighbourhood of the point $a$ in $\C$, then
we set $m_a(\varphi)=0$ if $\varphi(a)\neq 0$,
$m_a(\varphi)=+\infty$ if $\varphi=0$ in some neighbourhood of
$a$, and let $m_a(\varphi)$ be the multiplicity of $a$ if it is an
isolated zero of $\varphi$.

\begin{lemma}\label{var2} Let $x\in X$, ${\alpha} \in (-\infty,0)$
and  assume that $g\in \O(\overline \D,X)$, $g(0)=x$, and
$G_{g^*A}(0)<\alpha$.  Then there exist $f\in \O(\overline \D,X)$
and finitely many different points $a_1,\dots,a_k\in
\D\setminus\set 0$ such that $f(0)=x$ and
\begin{equation}
-\infty <\sum_{j=1}^k \tilde\nu_{f^*A}(a_j)\log|a_j|<\alpha.
\label{eq:var0}
\end{equation}
\end{lemma}

\begin{proof}  We have $G_{g^*A}(0)=\sum_{a\in
\D}\tilde\nu_{g^*A}(a)\log|a|<\alpha$, so we can choose finitely
many points $a_1,\dots,a_k\in \D\setminus\set 0$ such that
\begin{equation}
\sum_{j=1}^k \tilde\nu_{g^*A}(a_j)\log|a_j|<\alpha.
\label{eq:var1}
\end{equation}
If the sum in (\ref{eq:var1}) is finite, we take $f=g$.  If the
sum is equal to $-\infty$  and $g(\D)$ is not contained in $|A|$,
then $a_j=0$ and $0<\tilde\nu_{g^*A}(a_j)<+\infty$ for some $j$.
We choose $a \in \D\setminus\set{0}$ so close to $0$ that
$\log|a|<{\alpha}$ and $g$ is holomorphic in a neighbourhood of
the image of $h:\overline\D\to \C$, $h(z)=z(z-a)$.    If
$\psi_1,\dots,\psi_m$ are local generators for $\I_A$ near $x$,
then $m_a(\psi_j\circ g\circ h)=m_0(\psi_j\circ g)$  for all $j$,
which implies $\tilde\nu_{g^*A}(a)=\tilde\nu_{(g\circ h)^*A}(a)$.
If we set $f=g\circ h$, $k=1$, and $a_1=a$, then $f(0)=x$ and
(\ref{eq:var0}) holds.

If the sum in  (\ref{eq:var1}) equals $-\infty$  and $g(\D)$ is
contained in $|A|$, then we may replace $g$ by the constant disc
$z\mapsto x=g(0)$. We choose a neighbourhood $U$  of $x$ in $X$
and a biholomorphic map $\Phi:U\to\D^n$ such that  $\Phi(x)=0$. We
take $v\in \C^n$ with $|v|<1$ such that the disc $\overline \D\to
X$, $z\mapsto \Phi^{-1}(zv)$ is not contained in $|A|$ and choose
$a\in \D\setminus\set{0}$ so small that $\log |a|<{\alpha}$ and
$z(z-a)v\in \D^n$ for all $z\in \overline \D$.  If we take $k=1$,
$a_1=a$, and let $f$ be the map $z\mapsto \Phi^{-1}(z(z-a)v)$,
then $f(0)=f(a)=x\in |A|$, $0< \nu_{f^*A}(a)<+\infty$, and
(\ref{eq:var0}) holds.
\end{proof}

\begin{prooftx}{Proof of Theorem \ref{th:product}}
We need to prove that $G_A(x)\leq
\max\set{G_{A_1}(x_1),G_{A_2}(x_2)}$. Take $\alpha\in (-\infty,0)$
larger than the right hand side of this inequality. It is then
sufficient to show that $G_A(x)<{\alpha}$.

By Theorem~\ref{th:envelope} and  Lemma \ref{var2} we have $f_j\in
\O(\overline\D,X_j)$ with $f_j(0)=x_j$ and $a_{jk}\in
\D\setminus\set{0}$, $k=1,\dots,l_j$, $j=1,2$, such that
\begin{equation}
-\infty<\sum_{k=1}^{l_j}\tilde\nu_{f_j^*A_j}(a_{jk})\log|a_{jk}|<\alpha,\qquad
j=1,2. \label{eq:9.1}
\end{equation}
We choose $f_j$ so that $l_j$ becomes as small as possible. Then
$0<\tilde\nu_{f_j^*A_j}(a_{jk})<+\infty$  and $a_{jk}\neq 0$ for
all $j$ and $k$.   We define the Blaschke products $B_j$ by
$$
B_j(z)=\prod_{k=1}^{l_j}\bigg(\dfrac{a_{jk}-z}{1-\bar
a_{jk}z}\bigg)^{\mu_{jk}},  \qquad \text{ where } \
\mu_{jk}=\nu_{f_j^*A_j}(a_{jk}).
$$
Then  (\ref{eq:9.1}) implies $|B_j(0)|<e^\alpha$. We set
$b_j=B_j(0)$ and $\mu_j=\sum_{k=1}^{l_j}\mu_{jk}$ and we may
assume that $|b_1|\geq |b_2|$. We have
$B_j'(0)=B_j(0)\sum_{k=1}^{l_j}\mu_{jk}(|a_{jk}|^2-1)/a_{jk}$. If
$B_1'(0)=0$ we precompose $f_1$ with  a map $\D\to \D$ which fixes
the origin and makes a slight change of the points $a_{1k}$ so
that $B_1'(0)\neq 0$. By Schwarz Lemma this operation increases
the value of $|b_1|$, so we still have $|b_1|\geq|b_2|$.  By
precomposing $f_1$ by a rotation, we may assume that $B_1(0)=b_1$
is not a critical value of $B_1$.

If $c_j$ is one of the points $a_{jk}$ having largest absolute
value, then $|c_j|e^\beta\leq |b_j|$. For proving this inequality
we assume the reverse inequality $|b_j|<|c_j|e^\beta$  and for
simplicity enumerate the points so that
$|a_{j1}|\leq|a_{j2}|\leq\cdots$. Then
$$
\prod_{k=1}^{m_j}\bigg|\dfrac{a_{jk}}{c_j}\bigg|^{{\mu}_{jk}}<e^\beta
$$
where $m_j<l_j$ is the smallest natural number with
$|a_{jk}|=|c_j|$ for $k>m_j$.  Hence (\ref{eq:9.1}) holds with
$f_j$ replaced by $z\mapsto f_j(c_jz)$, $a_{jk}$ replaced by
$a_{jk}/c_j$, and $l_j$ by $m_j$, which contradicts the fact that
$l_j$ is minimal.

We may assume that $b_1=b_2$. Indeed, if $|b_1|>|b_2|$, we choose
$t\in (0,1)$ with $t^{-{\mu}_2}|b_2|=|b_1|$. Then $|a_{2k}|<t$,
for
$$
|a_{2k}|^{{\mu}_2}\leq |c_2|^{{\mu}_2}\leq |b_2|e^{-\beta}
<|b_2/b_1|=t^{{\mu}_2}.
$$
Replacing $f_2$ by $z\mapsto f_2(tz)$ and $a_{2k}$ by  $a_{2k}/t$,
we get $|b_1|=|b_2|$. Finally, replacing $f_2$ by $z\mapsto
f_2(e^{i\theta}z)$, where $e^{i\theta{\mu}_2}=b_2/b_1$, and
replacing $a_{2k}$  by $e^{-i\theta}a_{2k}$, we get $b_1=b_2$.

We let $C$ denote the set of all critical values of $B_1$.   We
have $B_1(0)=B_2(0)$, so we can take $\varphi_2:\D\to \D\setminus
B_2^{-1}(C)$ as the universal covering map with $\varphi(0)=0$. A
theorem of Frostman, see \cite{Nos}, p.~33, states that an inner
function on $\D$ omitting $0$ as a non-tangential boundary value
is a Blaschke product.  It is easy to show, see \cite{LarLasSig1},
p.~272, that since $0\not\in B_2^{-1}(C)$, $\varphi_2$ satisfies
the assumption in Frostman's theorem and is thus a Blaschke
product. The restriction of $B_1$ to $\D\setminus B_1^{-1}(C)$ is
a finite covering over $\D\setminus C$, so by lifting $B_2\circ
\varphi_2$ we conclude that there exists a function
$\varphi_1:\D\to B_1^{-1}(C)$ with $\varphi_1(0)=0$ and
$B_1\circ\varphi_1=B_2\circ\varphi_2$ and Frostman's theorem
implies again that $\varphi_1$ is a Blaschke product. Since
$|B_j\circ\varphi_j|=1$ almost everywhere on $\T$ and
$B_j(0)=b_j$, we can choose $r\in (0,1)$ such that
$$
\log|B_j\circ\varphi_j(0)|-\dfrac 1{2\pi}\int_0^{2\pi}
\log|B_j\circ\varphi_j(re^{i\theta})|\, d\theta <\alpha.
$$
We set $\sigma(z)=B_1\circ \varphi_1(rz)=B_2\circ\varphi_2(rz)$.
By the Poisson--Jensen representation formula, the left hands side
of this inequality equals $\sum_{i=1}^n\nu_i\log|z_i|$, where
$z_i$ are the zeros of $\sigma$ in $\D$ with multiplicities
$\nu_i$ for $i=1,\dots,n$.

We define $g_j\in \O(\overline \D,X_j)$ by
$g_j(z)=f_j\circ\varphi_j(rz)$ and $f\in \O(\overline \D,X)$ with
$f(0)=(x_1,x_2)$ by $f=(g_1,g_2)$. If $\sigma(z_i)=0$, then
${\varphi}_j(rz_i)=a_{jk_j}$ for some $k_j$, and
$$
\nu_i=m_{z_i}(\sigma)
=\mu_{jk_j}m_{z_i}(a_{jk_j}-\varphi_j(r\cdot))
=\tilde\nu_{f_j^*A_j}(a_{j,k_j})m_{z_i}(a_{jk_j}-\varphi_j(r\cdot))
=\tilde\nu_{g_j^*{A_j}}(z_i)
$$
Since the left hand side of this equation is independent of $j$,
we get
$$
\tilde\nu_{f^*A}(z_i)=\min_j\set{\tilde\nu_{g_j^*{A_j}}(z_i)}
=\nu_i.
$$
Hence
$$
G_A(x)\leq G_{f^*A}(0)=\sum_{a\in \D}\tilde\nu_{f^*A}(a)\log|a|
\leq \sum_{i=1}^n\tilde\nu_{f^*A}(z_i)\log|z_i|
=\sum_{i=1}^n\nu_i\log|z_i|<\alpha.
$$
\end{prooftx}

\section{Examples}

\begin{example}  Let $X$ be the unit polydisc $\D^n$ in $\Cn$,
$1\leq p\leq n$, and let $A$ be generated by
$\psi_k(z)=z_k^{\nu_k}$ for $1\leq k\leq p$ and positive integers
$\nu_k$.  Then the product property gives
$$
G_A(z)=\max_{1\leq k\leq p}\nu_k\log|z_k|.
$$
Furthermore, we have
$$
\big(dd^cG_A\big)^p=\nu_1\cdots\nu_p\, [|A|].
$$
\end{example}

\begin{example}
Let $X=\D^n$, $n\ge 2$, and let $A$ be generated by
$\psi_1(z)=z_1^2$, $\psi_2(z)=z_1z_2$. Then
$$G_A(z)=v(z):=\log|z_1|+\max\set{\log|z_1|,\log|z_2|}, \qquad
z=(z_1,z_2,z'')\in \D^n.
$$

First we take any $z\in \D^n\setminus\set{z_1=0}$ with $|z_1|\ge
|z_2|$ and consider the disc
$$ f(\zeta)=(\zeta,\frac{z_2}{z_1}\zeta,z''),\quad \zeta\in\D.
$$
Then for any $u\in\F_A$ we have
$f^*u(\zeta)\le\log|f^*\psi(\zeta)| +C=2\log|\zeta|+C$ and so,
since $u\le 0$, $f^*u(\zeta)\le 2\log|\zeta| =f^*v(\zeta)$. As
$f(z_1)=z$, this gives us $u(z)\le v(z)$.

For $z\in \D^n\setminus\set{z_1=0}$ with $|z_1|< |z_2|$, we take
the disc
$$ g(\zeta)=(\frac{z_1}{z_2}\zeta,\zeta,z''),\quad \zeta\in\D.
$$
Then for any $u\in\F_A$ we have again $g^*u(\zeta)\le
2\log|\zeta|+C$ near the origin and, since $u(z)\le\log|z_1|$
(which is the Green function for the polydisc with the poles along
the space $z_1=0$), $g^*u(\zeta)\le \log|z_1/z_2|$ near
$\partial\D$. Therefore, $g^*u(\zeta)\le
2\log|\zeta|+\log|z_1/z_2|=g^*v(\zeta)$ everywhere in $\D$. Since
$g(z_2)=z$, this shows $u(z)\le v(z)$ at all such $z$ as well.

Note that $f$ is, up to a M\"obius transformation, an extremal
disc for the disc functional $f\mapsto G_{f^*A}(0)$, while $g$ is
not. Note also that we have $dd^cG_A=[z_1=0]+Q$, where the current
$Q=dd^c\max\set{\log|z_1|,\log|z_2|}$ has the property
$Q^2=[z_1=z_2=0]$.
\end{example}

\begin{example}
Consider the variety $|A|=\{z_1=z_2=0\}\cup
\{z_2=z_3=0\}\cup\{z_1=z_3=0\} $ in the unit polydisk $\D^3$ of
$\C^3$. It is easy to see that the corresponding reduced complex
space $A$ is generated by $\psi(z)=(z_1z_2,z_2z_3,z_1z_3)$ and
that $|A|$ has the decomposition (in the sense of
Prop.~\ref{decomp-theo}) $|A|=J^2\cup J^3$ with $J^3=\set{0}$. We
claim that
$$G_A(z)=v(z):=\max\,\{\log|z_1z_2|,\log|z_2z_3|,\log|z_1z_3|\}.
$$
It suffices to check the relation $u(z)\le v(z)$ for any function
$u\in\F_A$ and each point $z\in\D^3$ with $|z_1|\ge |z_2|\ge
|z_3|$, $z_2\neq 0$. We take first any $z$ with $|z_1|= |z_2|\ge
|z_3|$ and consider the disc $f(\zeta)=\zeta z/|z_1|$,
$\zeta\in\D$. Then $f^*u\in SH^-(\D)$ and, since
$u\le\log|\psi|+C_1$ near the origin, $f^*u(\zeta)\le
\log|f^*\psi(\zeta)|+C_1=2\log|\zeta|+C_2$ when
$|\zeta|<\epsilon$. Therefore, $f^*u(\zeta)\le
2\log|\zeta|=f^*v(\zeta)$ and, in particular, $u(z)=f^*u(|z_1|)\le
f^*v(|z_1|)=v(z)$. The disc $f$ is, up to a M\"obius
transformation, an extremal disc for the disc functional $f\mapsto
G_{f^*A}(0)$ at such a point $z$.

Now we can take any $z$ with $|z_1|> |z_2|\ge |z_3|$, $|z_2|\neq
0$, and consider the analytic disc $g(\zeta)=(z_1, \zeta z_2,\zeta
z_3)$, $\zeta\in D_R$ with $R=|z_1|/|z_2|>1$. We have
$|g_1(\zeta)|= |g_2(\zeta)|\ge |g_3(\zeta)|$ when $|\zeta|=R$ and
thus $g^*u\le g^*v$ on $\partial D_R$. Furthermore,
$g^*u(\zeta)\le \log|g^*\psi(\zeta)|+C_3\le \log|\zeta|+C_4$ near
the origin. Since $g^*v(\zeta)=\log|\zeta z_1z_2|$, this shows
that $g^*u\le g^*v$ on $D_R$. Hence we get $u(z)=g^*u(1)\le
f^*v(1)=v(z)$, which proves the claim.

The current $(dd^c G_A)^2$ has Lelong numbers equal $1$ at each
point $a\in J^2=|A|\setminus\{0\}$. The point $0$ is exceptional:
the Lelong number $\nu(dd^cG_A,0)=2$, so $\nu((dd^cG_A)^2,0)\ge
4$, while $\nu([|A|],0)=3$.
\end{example}

\begin{example}
Let $X$ be the unit ball $\B_n$ in $\Cn$, $1\leq p\leq n$, and let
$A$ be generated by $\psi_k(z)=z_k$ for all $1\leq k\leq p$. In
the notation $z=(z',z'')$ with $z'\in\C^p$ and $z''\in\C^{n-p}$,
the Green function
$$
G_A(z)=\log\frac{|z'|}{\sqrt{1-|z''|^2}},
$$
because its restriction to every plane $z''=c\in\B_{n-p}$ is the
pluricomplex Green function for the ball of radius
$\sqrt{1-|c|^2}$ in $\C^p$ with simple pole at the origin.
\end{example}

\begin{example} The Green function $G_A$ for the unit ball $\B_n$
in $\Cn$, $n\ge 2$, with respect to $A$ generated by
$(\psi_1,\psi_2)=(z_1^2,z_2)$ is given by
$$G_A(z)=\dfrac 1 2\log\bigg( \frac{|z_1|^4}{(1-|z''|^2)^2}
+ \frac{2|z_2|^2}{1-|z''|^2} + \frac{|z_1|^2}{1-|z''|^2} \sqrt{
\frac{|z_1|^4}{(1-|z''|^2)^2} + \frac{4|z_2|^2}{1-|z''|^2} }\bigg)
 - \dfrac 12\log 2,$$
for $z=(z_1,z_2,z'')\in \B_n$ (compare with the formula for the
pluricomplex Green function with two poles in the ball
\cite{Coman}).

For proving this we let $v$ denote the function defined by the
right hand side.  Then $v\in \PSH(\B_n)\cap C(\overline\B_n
\setminus |A|)$ and satisfies $v(z)\le
\max\{\log|z_1|^2,\log|z_2|\}+C$ locally near
$|A|=\set{z_1=z_2=0}$.

Let us show that its boundary values on $\partial \B_n\setminus
|A|$ are zero. Take any $z\in\partial \B_n\setminus |A|$, then
$|z_1|^2=a$, $|z_2|^2=b$, $|z''|^2=1-a-b$ with $a,b\ge 0$,
$0<a+b\le 1$. We get
$$
v(z)= \dfrac 12\log\left[ \frac{a^2}{(a+b)^2} + \frac{2b}{a+b} +
\frac{a}{a+b}\left(2-\frac{a}{a+b}\right) \right]
 - \dfrac 12\log 2=0.
$$

 Finally we show that $v(z)\ge G_A(z)$ for almost all $z\in \B_n$
 (which implies $v\equiv G_A$). Take any $z\in \B_n$ with $z_1\neq 0$
 and consider the analytic curve
 $$f(\zeta)=(\zeta, \frac{z_2}{z_1^2}\zeta^2, z'').
 $$
Note that $f(z_1)=z$. We have $f^*v(\zeta) =2\log(|\zeta|/R(z))$,
while $f^*G_A$ is a negative subharmonic function in the disc
$|\zeta|<R(z)$ with the singularity $2\log|\zeta|$. So $f^*G_A \le
f^*v$ and, in particular, $G_A(z)=f^*G_A(z_1)\le f^*v (z_1)=v(z)$.

This shows also that $f$ is, up to a M\"obius transformation, an
extremal disc for the disc functional $f\mapsto G_{f^*A}(0)$ at
$z$ with $z_1\neq 0$. A corresponding extremal curve for
$z=(0,z_2,z'')$ is $f(\zeta)= (0,\zeta,z'')$.
\end{example}

\medskip
{\small{\it Acknowledgments.} The authors thank Daniel Barlet and
Alain Yger for valuable discussions. Part of the work was done
during Alexander's visit to the University of Iceland and Ragnar's
visit to H\o gskolen i Stavanger, and the authors thank the both
institutions for their kind hospitality.}

\bigskip
Alexander Rashkovskii

Tek/nat, H\o gskolen i Stavanger, POB 8002, 4068 Stavanger, Norway

E-mail: alexander.rashkovskii@his.no

\medskip

Ragnar Sigurdsson

Science Institute, University of Iceland, Dunhaga 3, IS-107
Reykjavik, Iceland

E-mail: ragnar@hi.is

\end{document}